\documentclass{amsart}

\usepackage{amsmath,amsfonts,amssymb,amsthm,amscd,latexsym,euscript}
\usepackage[all]{xy}

\newcommand{\dgrm}[1]{\ensuremath{\smash{\underset{\widetilde{\hphantom{#1}}}{#1}} \mathstrut}}

\DeclareMathOperator{\domain}{\textup{dom}}
\DeclareMathOperator{\codom}{\textup{codom}}

\newcommand{\cal}[1]{\ensuremath{\mathcal #1}}

\newtheorem {theorem1}{Theorem}[section]
\newtheorem {theorem}[theorem1]{Theorem}
\newtheorem {corollary}[theorem1]{Corollary}
\newtheorem {proposition}[theorem1]{Proposition}

\theoremstyle{definition}
\newtheorem {definition}[theorem1]{Definition}

\theoremstyle{remark}
\newtheorem {remark}[theorem1]{Remark}

\newcommand{\cat}[1]{\ensuremath{\EuScript #1}}

\newcommand{\Ch}{\ensuremath{\textup{Ch}}}

\DeclareMathOperator{\Map}{\ensuremath{\textup{Map}}}
\DeclareMathOperator{\Hom}{\ensuremath{\textup{Hom}}}
\DeclareMathOperator{\map}{\textup{map}}

\newcommand{\NN}{\ensuremath{\mathbb{N}}}
\newcommand{\ZZ}{\ensuremath{\mathbb{Z}}}

\newcommand{\colim}{\ensuremath{\mathop{\textup{colim}}}}

\def\id{\ensuremath{\mathrm{Id}}}

\newcommand{\rarrow}{\rightarrow}

\newcommand{\etale}{\'etale }

\newcommand{\Adamek}{Ad\'amek }
\newcommand{\Rosicky}{Rosick\'y }

\newcommand{\Mardesic}{Marde{\v{s}}i{\'c} }

\begin{document}

%Use Computer Modern fonts for the arrow heads
\SelectTips{cm}{10}

\title  [Generalized small object argument]
        {A generalization of Quillen's small object argument}
\author{Boris Chorny}

\address{Department of Mathematics, Middlesex College, The University of Western Ontario, London, Ontario N6A 5B7, Canada}

\email{bchorny2@uwo.ca}



\date{\today}
\dedicatory{}
\commby{}

\begin{abstract}
We generalize the small object argument in order to allow for its application to proper classes of maps (as opposed to sets of maps in Quillen's small object argument). The necessity of such a generalization arose with appearance of several important examples of model categories which were proven to be non-cofibrantly generated \cite{AHRT, Chorny, ChristHov, Isaksen-sharp}. Our current approach allows for construction of functorial factorizations and localizations in the equivariant model structures on diagrams of spaces \cite{Farjoun} and diagrams of chain complexes.  We also formulate a non-functorial version of the argument, which applies in two different model structures on the category of pro-spaces \cite{EH76, Isaksen-sharp}.

The examples above suggest a natural extension of the framework of cofibrantly generated model categories. We introduce the concept of a class-cofibrantly generated model category, which is a model category generated by classes of cofibrations and trivial cofibrations satisfying some reasonable assumptions.
\end{abstract}

\maketitle

\section{Introduction}
Quillen's definition of a model category has been slightly revised over the last decade. The changes applied to the first axiom MC1 requiring the existence of all finite limits and colimits, and to the last axiom MC5 requiring the existence of factorizations. The modern approaches to the subject \cite{Hirschhorn, Hovey} demand the existence of all small limits and colimits in MC1. This gives some technical advantages while treating transfinite constructions, such as localizations, in model categories. The modern version of the axiom MC5 requires the factorizations to be functorial. Unfortunately we will not be able to accept this stronger form of MC5, since we will be interested also in model categories for which the existence of functorial factorizations is unknown.

The most widely known model category without functorial factorizations is the category of pro-spaces or, more generally, of pro-objects (in the sense of Gro\-then\-dieck) in a proper model category \cat C \cite{EH76, Isaksen-strict} and its Bousfield localization modelling the \etale homotopy theory \cite{Artin-Mazur, Isaksen-sharp}. We introduce a new construction of factorizations in these model categories, but our construction still lacks functoriality.

The main tool for the construction of (functorial) factorizations in model categories and localizations thereof is Quillen's small object argument \cite{Hirschhorn, Hovey, Quillen}. However, in its original form, the argument is applicable neither to the category of diagrams with the equivariant model structure \cite{Farjoun}, nor to pro-categories, since it allows for the application in cofibrantly generated model categories only. We propose here a generalization which may be used in a wider class of model categories. The collections of generating cofibrations and generating trivial cofibrations may now form proper classes, satisfying the conditions of the following theorem:

\begin{theorem}[The generalized small object argument]\label{argument}
Suppose \cat C is a category containing all small colimits, and $I$ is a class of maps in \cat C satisfying the following conditions:
\begin{enumerate}
\item There exists a cardinal $\kappa$, such that each element $A\in \domain(I)$ is $\kappa$-small relative to $I$-cof;
\item
For every map $f\in \Map \cat C$ there exists a (functorially assigned) map $S(f)\in I\text{\textrm{-cof}}$ equipped with a (natural) morphism of maps $t_f\colon S(f)\rarrow f$, such that any morphism of maps $i\rarrow f$ with $i \in I$ factors through the (natural) map $t_f$.
\end{enumerate}
Then there is a (functorial) factorization $(\gamma,\delta)$ on \cat C such that, for all morphisms $f$ in \cat C, the map $\gamma(f)$  is in $I$-cof and the map $\delta(f)$ is in $I$-inj.
\end{theorem}

\begin{remark}
Note that the theorem above contains two sets of conditions and two statements: about existence of functorial and non-functorial factorizations. In the rest of the paper we refer to these statements as to the functorial and the non-functorial versions of the argument respectively. The proof of the theorem is given only for the functorial version. The adaptation of the proof for the non-functorial version may be achieved by removing the verification of the functoriality.
\end{remark} 

We say that a class $I$ of maps in a category \cat C \emph{permits the generalized small object argument} if it satisfies conditions (1) and (2) of Theorem~\ref{argument}.

This theorem is the second attempt by the author to generalize the small object argument. The previous version appeared in the study of the equivariant localizations of diagrams of spaces \cite{PhDI}. The specific properties of the equivariant model category of $D$-shaped diagrams of spaces and also the non-functorial factorization technique developed by E.~Dror Farjoun in \cite{Farjoun} suggested a rather complicated technical notion of instrumentation. It is essentially a straightforward ``functorialization'' of Dror Farjoun's ideas, which contain, implicitly, the non-functorial version of the generalized small object argument. The classes of generating cofibrations and generating trivial cofibrations of diagrams satisfy the conditions of instrumentation, but it is difficult (or impossible) to verify these conditions in other model categories. The conditions of Theorem~\ref{argument} on the class $I$ of maps are easier to handle and also more general then those of instrumentation, as we explain in Section~\ref{diagrams}. 

This paper shows also that two rather different homotopy theories of pro-spaces and of diagrams of spaces fit into a certain joint framework. In order to describe the similarity between the them let us give the following 
\begin{definition}
A model category \cat C is called \emph{class-cofibrantly generated} if 
\begin{enumerate}
\item there exists a class $I$ of maps in \cat C (called a class of \emph{generating cofibrations}) that permits (either functorial or non-functorial version of) the generalized small object argument and such that a map is a trivial fibration if and only if it has the right lifting property with respect to every element of $I$, and
\item there exists a class $J$ of maps in \cat C (called a class of \emph{generating trivial cofibrations}) that permits (either functorial or non-functorial version of) the generalized small object argument and such that a map is a fibration if and only if it has the right lifting property with respect to every element of $J$.
\end{enumerate}
\end{definition} 

The categorical dual to a class-cofibrantly generated model category is called \emph{class-fibrantly generated}.

The purpose of this paper is to give several non-trivial examples of class-cofibrantly and class-fibrantly generated model categories. One of the model categories we discuss here is new, others are classical and thus we only give a new construction of factorizations, applying the current version of the small object argument.

In particular we show that the equivariant model structure on the diagrams of spaces is class-cofibrantly generated, construct the equivariant model structure on the diagrams of chain complexes and prove that both known model structures on the category of pro-spaces are class-fibrantly generated.

Although the non-functorial version of the argument proves a weaker result, it is useful in many model categories, where the existence of functorial factorizations is still an open question. Recent applications of the non-functorial version of the generalized small object argument include new model structures on pro-spaces \cite{Isaksen-completions}, pro-spectra \cite{Isaksen-pro-spectra} and on small diagrams over a large category \cite{Chorny-smalldiag}.

The applications of Quillen's small object argument are not limited to abstract homotopy theory. A similar argument is used, for example, in the theory of categories to construct reflections in a locally presentable category with respect to a small orthogonality class \cite[1.36]{AR}. Recently another generalization of the small object argument was considered by the category theorists J.~\Adamek, H.~Herrlich, J.~\Rosicky and W.~Tholen \cite{AHRT-small-object}. Their version of the argument applies to the ``injective subcategory problem'' in locally ranked categories -- a generalization of the notion of a locally presentable category which includes topological spaces. We hope that our  generalization of the small object argument will be applicable to the ``orthogonal subcategory problem" and ``injective subcategory problem'' with respect to some reasonable classes of morphisms.

\paragraph{\bf The rest of the paper is organized as follows:}
Section~\ref{Proof} is devoted to the proof of the generalized small object argument. The non-functorial version of the argument is left for the reader. Next, we review some of our previous results about the diagrams of spaces in Section~\ref{diagrams} and show how they fit into the newly established framework. We extend this approach to the diagrams of chain complexes in Section ~\ref{chain-compl}. After providing the necessary preliminaries on pro-categories in Section~\ref{prelim} we apply the generalized cosmall object argument in Section~\ref{applications}.

\medskip
\paragraph{\emph{Acknowledgements.}} I would like to thank Dan Isaksen for many fruitful conversations and suggestions for improving this paper. In particular I owe him the idea of the proof of Theorem~\ref{second}.
 
\section{Proof of the generalized small object argument}\label{Proof}
\begin{proof}[Proof of Theorem~\ref{argument}]
Given a cardinal $\kappa$ such that every domain of $I$ is $\kappa$-small
relative to $I$-cof, we let $\lambda$ be a $\kappa$-filtered ordinal (i.e., $\lambda$ is a limit ordinal and, if $\Lambda\subset \lambda$ and $|\Lambda|\le \kappa$, then $\mathop{\mathrm{sup}}\Lambda < \lambda$).
\par
To any map $f\colon X \rarrow Y$ we will associate a functor $Z^f\colon
\lambda \rarrow \cat C$ such that $Z_0^f = X$, and a natural transformation
$\rho^f\colon Z^f \rightarrow Y$ factoring $f$, i.e., for each $\beta <
\lambda$ the triangle
\[
\xymatrix@ur{
                                       X\ar[dr]^f \ar[d]\\
Z^f_\beta \ar[r]_{\rho^f_\beta}      & Y }
\]
is commutative. Each map $i^f_\beta\colon  Z_\beta^f \rightarrow
Z_{\beta+1}^f$ will be a pushout of a map of the form $S(f)$, i.e., $i^f_\beta
\in I$-cof, since $I$-cof is closed under pushouts.
\par
We will define $Z^f$ and $\rho^f\colon Z^f \rarrow Y$ by transfinite induction,
beginning with $Z^f_0 = X$ and $\rho_0^f = f$. If we have defined
$Z_{\alpha}^f$ and $\rho_{\alpha}^f$ for all $\alpha < \beta$ for some limit
ordinal $\beta$, define $Z_{\beta}^f = \colim_{\alpha<\beta} Z_{\alpha}^f$,
and  define $\rho_\beta^f$ to be the map induced, naturally, by the
$\rho_\alpha^f$. Having defined $Z_\beta^f$ and $\rho_\beta^f$, we define
$Z_{\beta+1}^f$ and $\rho_{\beta+1}^f$ as follows. Consider the natural map
$t(\rho_\beta^f)\colon S(\rho_\beta^f)\rarrow \rho_\beta^f$, i.e. the
following commutative square:
\[
\begin{CD}
A  @>{\domain(t(\rho_\beta^f))}>> Z_\beta^f \\
@V{S(\rho_\beta^f)}VV  @V{\rho_\beta^f}VV\\
B  @>>{\codom(t(\rho_\beta^f))}> Y.
\end{CD}
\]
Define $Z_{\beta+1}^f$ to be the pushout of this diagram and define
$\rho_{\beta+1}^f$ to be the map naturally induced by $\rho_{\beta}^f$.
\par
For each morphism $g=(g^1,g^2)\colon f_1\rarrow f_2$ in the category $\Map \cat
C$, i.e., for each commutative square
\[
\begin{CD}
{X_1}      @>{g^1}>>    X_2 \\
@V{f_1}VV            @V{f_2}VV\\
{Y_1}      @>{g^2}>>    {Y_2}
\end{CD}
\]
we define a natural transformation $\xi^g\colon Z^{f_1}\rarrow Z^{f_2}$ by
transfinite induction over small ordinals, beginning with $\xi^g_0 = g^1$. If
we have defined $\xi^g_\alpha$ for all $\alpha < \beta$ for some limit ordinal
$\beta$, define $\xi^g_\beta = \colim_{\alpha < \beta} \xi^g_\alpha$. Having
defined $\xi^g_\beta$, we define $\xi^g_{\beta+1}\colon Z^{f_1}_{\beta+1}
\rarrow Z^{f_2}_{\beta+1}$ to be the \emph{natural} map induced by $g_\beta =
(\xi^g_\beta, g^2)\colon  \rho^{f_1}_\beta \rarrow \rho^{f_2}_\beta$, namely
the \emph{unique} map between the pushouts of the horizontal lines of the
following diagram which preserves its commutativity:
\[
\begin{CD}
{B_1}  @<{S(\rho_\beta^{f_1})}<<  {A_1}    @>{\domain(t(\rho_\beta^{f_1}))}>>    Z_\beta^{f_1} \\
@V{h_2}VV       @V{h_1}VV    @VV{\xi_\beta^g}V\\
{B_2}  @<<{S(\rho_\beta^{f_2})}<  {A_2}    @>>{\domain(t(\rho_\beta^{f_2}))}>    Z_\beta^{f_2}. \\
\end{CD}
\]

In this diagram $(h_1, h_2)=S(g_\beta)$. The commutativity of the diagram
follows readily, since $S$ is a functor and $t$ is a natural transformation.
\par
The required functorial factorization $(\gamma,\delta)$ is obtained when we
reach the limit ordinal $\lambda$ in the course of our induction. Then we
define $\gamma(f)\colon  X \rarrow Z^f_{\lambda}$ to be the (transfinite)
composition of the pushouts, and $\delta(f) = \rho^f_{\lambda}\colon
Z^f_{\lambda} \rarrow Y$. $\gamma(f)\in I$-cof since $I$-cof is closed under
transfinite compositions.
\par
To complete the definition of the functorial factorization (see
\cite[1.1.1]{Hovey}, \cite[1.1.1]{Hovey_err}) we need to define for each
morphism $g\colon  f_1 \rarrow f_2$ a natural map $(\gamma,\delta)^g\colon
Z^{f_1}_{\lambda} \rarrow Z^{f_2}_{\lambda}$ which makes the appropriate
diagram commutative. Take $(\gamma,\delta)^g = \xi^g_\lambda$.
\par
It remains to show that $\delta(f) = \rho^f_{\lambda}$ has the right lifting
property with respect to $I$. To see this, suppose we have a commutative
square as follows:
\[
\begin{CD}
  C             @>h'>>              Z^f_\lambda \\
@VlVV                        @VV{\rho^f_\lambda}V\\
  D             @>k'>>              Y
\end{CD}
\]
where $l$ is a map of $I$. Due to the first condition of the theorem the object
$C$ is $\kappa$-small relative to $I$-cof, i.e., there is an ordinal $\beta <
\lambda$ such that $h'$ is the composite $C\stackrel{h_\beta}
{\longrightarrow} Z_\beta^f \longrightarrow Z^f_\lambda$. Hence we obtain the
following commutative diagram:
\[
\xymatrix{
C \ar[rr]^{h_\beta} \ar[dd]_l  &      &      Z^f_\beta \ar[d] \ar@(r,r)[dd]^{\rho^f_\beta}\\
                               &      &      Z^f_\lambda \ar[d]_{\rho^f_\lambda}\\
D \ar[rr]_{k'}                 &      &      Y. }
\]
The second condition of the theorem implies that there exists a factorization
in the category $\Map \cat C$ of the map $(h_\beta, k')$ through
$t(\rho^f_\beta)$ which is a map of maps with domain $S(\rho^f_\beta)\colon A
\rarrow B$ and range $\rho^f_\beta$, i.e., there is a commutative diagram
\[
\xymatrix{
C \ar[r]\ar@(u,u)[rr]^{h_\beta} \ar[dd]_l  &   A\ar[dd]_{S(\rho^f_\beta)}\ar[r]^{h}   &      Z^f_\beta \ar[d] \ar@(r,r)[dd]^{\rho^f_\beta}\\
                                           &       &      Z^f_\lambda \ar[d]_{\rho^f_\lambda}\\
 D \ar[r]\ar@(d,d)[rr]_{k'}                 &   B \ar[r]^{k} &      Y }
\]
where $(h,k) = t(\rho^f_\beta)$.

By construction, there is a map $B \stackrel{k_\beta}{\longrightarrow}
Z_{\beta+1}^f$ such that $k_\beta S(\rho_\beta^f) = i_\beta^f h$ and $k =
\rho^f_{\beta+1}k_\beta$, where $i_\beta^f$ is the map $Z_{\beta}^f\rarrow
Z_{\beta+1}^f$. The composition $D \longrightarrow B \stackrel
{k_\beta}{\longrightarrow} Z_{\beta+1}^f\longrightarrow Z^f_\lambda$ is the
required lift in the initial commutative square.
\end{proof}

\begin{remark}
In all the applications we have in mind, the map $S(f)$ is a coproduct of maps from $I$. Hence the construction above provides us with the factorization of any map $f$ into an $I$-cellular map followed by an $I$-injective map, as in the classical construction. But we prefer to leave the formulation of conditions on the class $I$ of maps in the present (simpler) form, since we hope that they will be useful elsewhere and we do not see a big advantage in $I$-cellular maps instead of $I$-cofibrations.
\end{remark}

\section{Example: Instrumented classes of maps permit the generalized small object argument}\label{diagrams}
Let us recall, in an informal manner, the notion of instrumentation introduced in \cite{PhDI}. Instrumentation for a class $I$ of maps in a category \cat C is a formalization of the following functorial version of the classical cosolution-set condition: for any morphism $f$ in \cat C there is a {\em naturally} assigned set of maps $\cal I(f) = \{i \rarrow f \;|\; i\in I\}$, such that for any morphism of maps $j\rarrow f$ with $j\in I$ there exists a factorization $j\rarrow i\rarrow f$ with $(i\rarrow f)\in \cal I(f)$. Additionally, every domain of a map in $I$ is $\kappa$-small with respect to $I$-cell for some fixed cardinal $\kappa$. 
 
\begin{proposition}\label{generalization}
Any instrumented class of maps $I$ in a category \cat C permits the generalized small object argument. 
\end{proposition}
\begin{proof}
The first condition of Theorem~\ref{argument} is satisfied because of the same assumption for instrumented classes of maps.

Instrumentation gives rise to the augmented functor $S$ in the following way:
\[
S(f) = \coprod \domain(\cal I(f)) = \coprod \left\{i \;|\; (i \rarrow f)\in \cal I(f) \right\}.
\]
Naturality of $\cal I$ ensures the functoriality of $S$. The augmentation $t_f\colon S(f)\rarrow f$ exists since every $i$ is equipped with a map into $f$, hence their coproduct is naturally mapped into $f$. Certainly $S(f)\in I$-cof, and the factorization property follows from the similar property of instrumentation.
\end{proof}

Instrumented classes of maps were applied to the study of \emph{equivariant} model structures on diagrams of spaces. Let $D$ be a small category. It was essentially shown in \cite{PhDI} that the category of $D$-shaped diagrams of spaces (by the category of spaces we mean here either the category of simplicial sets, or the category of compactly generated topological spaces with the standard simplicial model structure) is class-cofibrantly generated. In this section we show that the classes of generating cofibrations and generating trivial cofibrations permit (the new version of) the generalized small-object argument.

Let us recall first the definition of the equivariant model structure on $\cal S^D$ initially introduced in \cite{Farjoun}. We use the word \emph{collection} to denote a set or a proper class with respect to some fixed universe $\mathfrak U$. A $D$-diagram $\dgrm O$ of spaces is called an \emph{orbit} if $\colim_D \dgrm O = \ast$. We denote by $\cal O_D$ the collection of all orbits of $D$ (which is not necessarily a set). For any diagram $\dgrm W$ and a map $f\colon \dgrm X \rarrow \dgrm Y$, there is an induced map of simplicial sets $\map(\dgrm W,f)\colon \map(\dgrm W, \dgrm X) \rarrow \map(\dgrm W,\dgrm Y)$; see \cite{DF} for details.
\begin{definition}
In the equivariant model structure on ${\cal S}^D$ a morphism $f\colon \dgrm X \rarrow \dgrm Y$ is a
\begin{itemize}
\item
\emph{weak equivalence} if and only if $\map(\dgrm O,f)$ is a weak equivalence of spaces for any orbit $\dgrm O$; 
\item
\emph{fibration} if and only if $\map(\dgrm O,f)$ is a fibration of spaces for any orbit $\dgrm O$;
\item
\emph{cofibration} if and only if it has the left lifting property with respect to any trivial fibration.
\end{itemize}
\end{definition}
The standard axioms of simplicial model categories were verified in \cite{Farjoun} for the equivariant model structure on $\cal S^D$. Functorial factorizations were constructed in \cite{PhDI}. In that construction we used a different version of the generalized small-object argument, which applied only to instrumented classes of maps. The purpose of this section is to prove Proposition~\ref{generalization},which shows that Theorem~\ref{argument} provides a more general version of the argument.

Let $I = \{\dgrm O\otimes \partial \Delta^n \hookrightarrow \dgrm O \otimes \Delta^n \;|\; {\dgrm O\in\cal O_D,\; n\geq 0}\}$ and $J = \{\dgrm O\otimes \Lambda^n_k \mathbin{\tilde{\hookrightarrow}} \dgrm O\otimes \Delta^n \;|\; {\dgrm O\in\cal O_D,\; n\geq k \geq 0}\}$ be two classes of maps in $\cal S^D$. If the index category $D$ is such that $\cal O_D$ is a set (this happens, for example, when $D$ is a group), then the collections $I$ and $J$ are {\em sets} of cofibrations and the equivariant model structure on $\cal S^D$ is cofibrantly generated with $I$ equal to the set of generating cofibrations and $J$ equal to  the set of generating trivial cofibrations. But usually $I$ and $J$ are proper classes of maps, and it was shown in \cite{PhDI} that they form classes of generating cofibrations and generating trivial cofibrations in the equivariant model structure on the $D$-shaped diagrams of spaces, which is class-cofibrantly generated. Our aim here is to show that the classes $I$ and $J$ permit the functorial version of the generalized small object argument. This follows from Proposition~\ref{generalization}, since the classes $I$ and $J$ are both instrumented. 

\section{Equivariant model structure on the diagrams of chain complexes}\label{chain-compl}
Let $R$ be a commutative ring. The category $\Ch(R)$ of (unbounded) chain complexes of $R$-modules carries a cofibrantly generated model structure with weak equivalences being quasi-isomorphisms and fibrations being levelwise surjections of chain complexes \cite{Hinich, Hovey}. The purpose of this section is to extend this model structure to the category of diagrams of chain complexes with \emph{equivariant} weak equivalences (as opposed to objectwise).

Let $\cal A$ be a complete and cocomplete abelian category and $\cal P$ be a \textbf{projective class} in $\cal A$; see \cite{ChristHov} for the definition. Then, under certain conditions on $\cal P$, the category of chain complexes in $\cal A$ carries the relative model structure: a map $f\colon A\rarrow B$ is a weak equivalence or a fibration if the induced map $\hom (P, f)\colon \hom(P, A)\rarrow \hom(P,B)$ is a weak equivalence or a fibration in the standard model structure on $\Ch(\ZZ)$. We exhibit bellow a projective class in the category $\cal A$ of $D$-shaped diagrams of $R$-modules which satisfies the technical conditions that ensure the existence of the relative model structure.

Let $D$ be a small category. A $D$-diagram of sets \dgrm T is called an orbit if, as before, $\colim_D \dgrm T = \ast$. And let $\cal O_D$ be the collection of all $D$-orbits. For every orbit \dgrm T we associate a diagram of free $R$-modules $P_{\dgrm T} = R(\dgrm T)$. Recall that the free functor is the left adjoint of the forgetful functor $U(\,\cdot\,)$.

Let $\cal P' = \{P_{\dgrm T} \,|\, \forall \, \dgrm T\in \cal O_D\}$, and let $\cal P$ be the projective class determined by $\cal P'$, i.e., if $\cal E$ is the class of $\cal P'$-epic maps, then $\cal P$ is precisely the class of all diagrams P of chain complexes such that each map in $\cal E$ is $P$-epic. In other words, $\cal P'$ is a collection of enough projectives for $\cal P$.

The same argument as in \cite[3.1]{PhDI} shows that $P_{\dgrm T}$ are $\aleph_0$-small relative to split monomorphisms with $\cal P$-projective cokernel. Therefore the projective class $\cal P$ has enough $\aleph_0$-small projectives.

In order to conclude that the $\cal P$-relative model structure on the category $\Ch(\cal A)$ exists it suffices, by \cite[2.2(B)]{ChristHov}, to prove that $\cal P$-resolutions can be chosen functorially. We suggest the following construction: for every diagram of $R$-modules $\dgrm X \in \cal A$, let $\cal O_{\dgrm X} = \{\ast \times_{U(\colim_D\dgrm X)} U\dgrm X \;|\; x\colon \ast \rarrow U(\colim_D\dgrm X)\}$ be the set of orbits over each point $x\in U(\colim_D \dgrm X)$. Every element $\dgrm T_x\in \cal O_{\dgrm X}$ is equipped with the projection map $\pi_x\colon \dgrm T_x \rarrow U\dgrm X$. Consider the collection of adjoint maps $\varphi_x \colon P_{\dgrm T_x} \rarrow \dgrm X$, for all $x\in U(\colim_D \dgrm X)$. We define the $\cal P$-resolution functor $P\dgrm X = \bigoplus_{x\colon \ast \rarrow U(\colim_D \dgrm X)} P_{\dgrm T_x}$ and the canonical $\cal P$-epic map $\varepsilon_{\dgrm X}\colon P\dgrm X \rarrow \dgrm X$ is induced by the maps $\varphi_x$. 

We have to verify that $\varepsilon_{\dgrm X}$ is indeed $\cal P$-epic. It suffice to show that $\varepsilon_{\dgrm X}$ is $\cal P'$-epic, i.e., for every map $P_{\dgrm T} \rarrow \dgrm X$ we must construct a factorization $P_{\dgrm T}\rarrow P\dgrm X \rarrow\dgrm X$. By adjointness this is equivalent to factorization of the map $\dgrm T\rarrow U\dgrm X$ through $UP\dgrm X$. The later factorization exists by construction of $P\dgrm X$: denote the induced map $\colim_D \dgrm T = \ast \rarrow U(\colim_D \dgrm X)$ by $x$, then there exists a factorization $\dgrm T\rarrow \dgrm T_x \rarrow U\dgrm X$ which can be further refined to $\dgrm T\rarrow \dgrm T_x\rarrow UP_{\dgrm T_x}\rarrow UP\dgrm X \rarrow U\dgrm X$, where the second map is the unit of the adjunction and the third map is the inclusion.

This finishes the proof that there exist the equivariant model structure on the category of diagrams of chain complexes. But our main motivation in this example is to find another class-cofibrantly generated model category. This model category is generated by classes $I = \{ \Sigma^{n-1}P_{\dgrm T} \rarrow D^n P_{\dgrm T} \,|\, n\in \ZZ, \, P_{\dgrm T}\in \cal P'\}$ and $J = \{ 0 \rarrow D^n P_{\dgrm T} \,|\, n\in \ZZ, \, P_{\dgrm T}\in \cal P'\}$ of generating cofibrations and generating trivial cofibrations respectively. Lemma \cite[5.5]{ChristHov} implies that $\cal P$-relative fibrations are those maps that have the right lifting property with respect to $J$, and $\cal P$-relative trivial fibrations are those maps that have the right lifting property with respect to $I$.

We only need to show that the classes $I$ and $J$ satisfy the functorial version of the generalized small object argument. The domains of the elements of $J$ are obviously $\aleph_0$-small. The domains of the elements of $I$ are $\aleph_0$-small by Lemma \cite[4.3]{ChristHov} since every element of $\cal P'$ is $\aleph_0$-small. It remains to verify the existence of an augmented functorial construction which associates to every map $f\colon X\rarrow Y$ between diagram of chain complexes (=chain complexes of diagrams of $R$-modules) a (trivial) cofibration $F(f)$ with a natural map $t\colon F(f)\rarrow f$, such that for every element of $i\in I$ ($j\in J$), any map $i\rarrow f$ factors through $t$. We will give the construction for $I$. The construction for $J$ is similar.
 
Given $i\colon \Sigma^{n-1}P_{\dgrm T}\rarrow D^n P_{\dgrm T}$, the maps from $i$ to $f$ are in bijective correspondence with commutative squares:
\[
\xymatrix{
P_{\dgrm T} \ar[r]\ar[d]  & \dgrm X_{n-1}\ar[d]^{f_{n-1}} \\
\dgrm Y_n \ar[r]^{d}          &  \dgrm Y_{n-1}. }
\]
Let $\dgrm W_n = \dgrm X_{n-1}\times_{\dgrm Y_{n-1}} \dgrm Y_n$, then the commutative squares above are in bijective correspondence with the maps $P_{\dgrm T} \rarrow \dgrm W_n$, which are in bijective correspondence with the maps $\dgrm T\rarrow U\dgrm W_n$, by adjointness. Let $\cal O_{\dgrm W_n}$ be the set of orbits over each point in $U\dgrm W_n$. Then we define 
\[
F(f) = \bigoplus_{n\in \ZZ}\left(\bigoplus_{\dgrm O\in \cal O_{\dgrm W_n}}\underset {D^n P_{\dgrm O}}{\overset {\Sigma^{n-1}P_{\dgrm O}} \downarrow}\right).
\]
The augmentation map is induced by the maps $P_{\dgrm O} \rarrow \dgrm W_n$ and the factorization property is readily verified. The functoriality of the constructions follows from the naturality on each step.

We have shown that the equivariant model structure on the category of diagrams of chain complexes is class-cofibrantly generated. One can show, by an argument similar to \cite{Chorny} that this model structure is not cofibrantly generated.

\section{Preliminaries on Pro-Categories}\label{prelim}
\begin{definition}
A small, non-empty category $I$ is {\bf cofiltering} if for every pair of
objects $i$ and $j$ there exists an object $k$ together with maps $k\rarrow i$
and $k\rarrow j$; and for every pair of morphisms $i\underset g {\overset f
\rightrightarrows} j$ there exists a map $h\colon k\rarrow i$ with $fh = gh$.
A diagram is said to be \textbf{cofiltering} if its indexing category is so.

For a category \cat C, the category {\bf pro-${\bf \EuScript C}$} has objects
all cofiltering diagrams in \cat C, and the set of morphisms from a pro-object
$X$ indexed by a cofiltering $J$ into a pro-object $Y$ indexed by a
cofiltering $I$ is given by the following formula:
\[
\Hom_{\text{pro-\cat C}}(X,Y)=\lim_i \colim_j \Hom_{\cat C}(X_j,Y_i).
\]
A {\bf pro-object} ({\bf pro-morphism}) is an object (morphism) of pro-\cat C.
The structure maps of diagrams which represent pro-objects are also called {\bf
bonding} maps and denoted by $b_{i_1,i_2}\colon Y_{i_1}\rarrow Y_{i_2}$.
\end{definition}

This definition of morphisms in a pro-category requires, perhaps, some
clarification. By definition, a pro-map $f\colon X \rarrow Y$ is a compatible
collection of maps $\{f_i\colon X\rarrow Y_i \;|\; i\in I,\; f_i\in \colim_j
\Hom_{\cat C}(X_j, Y_i) \}$. By construction of direct limits in the
category of sets, $\colim_j \Hom_{\cat C}(X_j, Y_i)$ is the set of equivalence
classes of the elements of $\coprod_j \Hom_{\cat C}(X_j, Y_i)$. If for each
equivalence class $f_i$ we choose a representative $f_{j,i}\colon X_j
\rarrow Y_i$, then we obtain a representative of the pro-morphism $f$.
This motivates the following alternative characterization of the morphisms in
pro-\cat C (cf.~\cite[2.1.2]{EH76}):

\begin{definition}
A {\bf representative} of a morphism from a pro-object $X$ indexed by a
cofiltering $I$ to a pro-object $Y$ indexed by a cofiltering $K$ is a function
$\theta\colon K\rarrow I$ (not necessarily order-preserving) and morphisms
$f_k\colon X_{\theta(k)}\rarrow Y_k$ in \cat C for each $k\in K$ such that if
$k\leq k'$ in $K$, then for some $i\in I$ with $i\geq \theta(k)$ and $i\geq
\theta(k')$ the following diagram commutes:
\[
\xymatrix{
 X_i\ar[r]^{b_{i,\theta(k')}}\ar[dr]_{b_{i,\theta(k)}}
        & X_{\theta(k')}\ar[r]^{f_{k'}} & Y_{k'} \ar[d]^{b_{k',k}}\\
        & X_{\theta(k)}\ar[r]^{f_k}   & Y_{k}.
  }
\]

A representative $(\phi,\{g_k\})$ {\bf rarefies} the representative $(\theta,\{f_k\})$ if for every $k\in K$, $\phi(k)\geq \theta(k)$ and there exists a bonding map $b_{\phi(k),{\theta(k)}}\colon X_{\phi(k)}\rarrow X_{\theta(k)}$ with $g_k = f_k\circ b_{\phi(k),{\theta(k)}}$.

A representative $(\theta, \{f_k\})$ is called {\bf strict} \cite[p.~36]{Friedlander} if $\theta$ is a functor and the maps $\{f_k\colon X_{\theta(k)}\rarrow Y_k \;|\; k\in K\}$ constitute a natural transformation $f\colon X\circ\theta \rarrow Y$. In other words, the representative of a morphism is strict if all $f_k$ fit into commutative squares of the following form:
\[
\xymatrix{
 X_{\theta(k')}\ar[r]^{f_{k'}} \ar[d]_{b_{\theta(k'),\theta(k)}} & Y_{k'} \ar[d]^{b_{k',k}}\\
 X_{\theta(k)}\ar[r]^{f_k}  & Y_{k}
  }
\]

A strict representative $(\theta, \{f_k\})$ is called {\bf levelwise} if the domain and range of $f$ are indexed by the same indexing category $I$ and $\theta = \id_I$. 
\end{definition}

\begin{remark}
Not every pro-map has a levelwise representative, but every pro-map may be reindexed, up to a pro-isomorphism into a pro-map equipped with levelwise representative (see, e.g., \cite[A.3.2]{Artin-Mazur}). Corollary~\ref{levelwise} below gives a brief proof of existence of functorial levelwise replacement. In \cite{Isaksen-levelwise} Dan Isaksen proves that the original construction by Artin-Mazur is functorial.
\end{remark}

The proof of the following standard proposition may be found, for example, in \cite[Ch.~1\S1]{MarSe82}.
\begin{proposition}
Two representatives $(\theta,\{f_k\})$ and $(\theta',\{f'_k\})$ are called \textbf{equivalent} if there exists a
representative $(\theta'',\{f''_k\})$ which rarefies both of them. This relation between representatives is an equivalence relation. The equivalence classes of representatives of morphisms from a pro-object $X$ into a pro-object $Y$ are in natural bijective correspondence with the elements of $\Hom_{\text{\rm pro-}\cat C}(X,Y)$.
\end{proposition}

An important technique in pro-categories is reindexing. We use two types of
reindexing: for maps and for objects. Their crucial property is functoriality.

The following theorem, proven in \cite{Meyer}, will provide us with a
functorial choice of a levelwise representative of a pro-morphism:
\begin{theorem}[C.~V. Meyer]\label{equiv}
\sloppy Let \cat C be any category, and let
\[
 F\colon \text{\rm pro-}(\Map \cat C) \longrightarrow \Map(\text{\rm pro-}\cat C)
\]
be the obvious functor. Then $F$ is fully faithful and essentially surjective,
i.e.,  the categories $\text{\rm pro-}(\Map \cat C)$ and $\Map(\text{\rm
pro-}\cat C)$ are equivalent.

\end{theorem}
Fix once and for all the functors which induce the equivalences:
\[
F\colon \text{\rm pro-}(\Map \cat C) \leftrightarrows \Map(\text{\rm pro-}\cat
C) : \! G.
\]
Beware that the pro-objects of \cite{Meyer} are indexed by cofiltered
categories that are not necessarily small. In this paper we consider only
small indexing categories. Nevertheless, Theorem~\ref{equiv} is still true, as explained in \cite[3.1, 3.5]{Isaksen-calc}.
\begin{corollary}\label{levelwise}
{
There exists a functor $L\colon \Map(\text{\rm pro-}\cat C)\rarrow \Map(\text{\rm pro-}\cat C)$ naturally isomorphic to the identity satisfying the following property. For every $f\in \Map(\text{\rm pro-\cat C})$ the domain and the range of $L(f)$ are indexed by the same indexing category $I$ and there exists a strict representative $(\theta, \{f_k\})$ of $f$ with $\theta = \id_I$. In other words, $f$ has a levelwise representative.\sloppy

}
\end{corollary}
\begin{proof}
Take $L = FG$.
\end{proof}

We are going to use inductive arguments, therefore we need a \emph{functorial} reindexing result that produces pro-objects indexed by \emph{cofinite strongly directed sets}

\begin{definition}
A partially ordered set $I$ is \textbf{directed} if for any $i,i'\in I$ there
exists $i''\in I$ with $i''\geq i$ and $i''\geq i'$; $I$ is called
\textbf{strongly} directed if $i\geq i'$ and $i\leq i'$ implies $i=i'$. A
directed set $I$ is called \textbf{cofinite} if every $i\in I$ has only a
finite number of predecessors.
\end{definition}

The following theorem \cite[2.1.6]{EH76}\cite[p.~15, Thm.~4]{MarSe82} supplies
us with the required reindexing for pro-objects:
\begin{theorem}\label{Mardesic}
There exists a functor $M\colon \text{\rm pro-}\cat C\rarrow \text{\rm
pro-}\cat C$ (called the \Mardesic functor), naturally equivalent to the identity,
such that $M(X)$ is indexed by a cofinite strongly directed set for every $X$
in {\rm pro-\cat C}.
\end{theorem}

The following proposition (see \cite[p.~9, Lemma~2]{MarSe82} for the proof)
allows us to assume that whenever our pro-objects are indexed by cofinite
strongly directed sets we may choose a strict representatives for each
morphism.

\begin{proposition}
Let $f\colon X\rarrow Y$ be a pro-map. If $X$ and $Y$ are indexed by directed
sets and the indexing category of $Y$ is cofinite and strongly directed, then
there exists a strict representative $(\theta, \{f_k\})$ of $f$.
\end{proposition}

\section{Applications of the generalized small object argument: Functorial factorizations in pro-$\cat C$}\label{applications}

We discuss in this section the application of Theorem~\ref{argument}, or more precisely of its dual, to the two different model categories on the category of pro-spaces. The strict model  structure was constructed by D.A.~Edwards and H.M.~Hastings \cite{EH76}. Weak equivalences and cofibrations are {\bf essentially levelwise} in the strict model structure, i.e., levelwise, up to a reindexing. The properness of \cat C is the only condition required for the existence of the strict model structure on pro-\cat C \cite{Isaksen-strict}. The localization of the strict model structure on the category of pro-(simplicial sets) with respect to the class of maps rendered into strict weak equivalences by the functor $P$ (the functor which replaces a space with its Postnikov tower) was constructed by D.~Isaksen \cite{Isaksen-sharp}. Weak equivalences in Isaksen's model structure generalize those of Artin-Mazur \cite{Artin-Mazur} and J.W.~Grossman \cite{Grossman}. 

We recall from \cite[2.1.7]{Hovey} that given a class $M$ of maps in a model category, $M$-proj is the class of maps which have the left lifting property with respect to every map in $M$.
 
\begin{theorem}\label{first}
Let \cat C be a proper model category. Then the strict model structure on {\rm pro-}\cat C is class-fibrantly generated by the classes $M$ of constant fibrations and $N$ of constant trivial fibrations. The classes of maps $M$ and $N$ admit the non-functorial version of the generalized cosmall object argument.
\end{theorem}
\begin{proof}
The class of trivial cofibrations equals $M$-proj and the class of cofibrations between pro-objects equals $N$-proj \cite[5.5]{Isaksen-strict}. Therefore we only have to prove that the classes $M$ and $N$ satisfy the non-functorial version of the  generalized cosmall object argument. It was shown in \cite{Isaksen-calc} that every constant pro-object is countably cosmall.

In order to apply the generalized cosmall object argument on the class $M$, we need to construct a map $S\colon \Map\cat C\rarrow \Map\cat C$ equipped with a morphism of maps $t\colon \mathrm{Id}_{\Map \cat C}\rarrow S$ such that for every $f\in \Map\cat C$, $S(f)$ is an $M$-fibration and any morphism of maps $f\rarrow g$ with $g\in M$ factors through $t(f)\colon f\rarrow S(f)$.

We start by replacing $f$ with an isomorphic levelwise representative $\{f_i\}_{i\in I}$ for some cofiltering $I$; Corollary~\ref{levelwise} allows to do it in a functorial way. Next, for every $i\in I$, we factorize $f_i$ into a trivial cofibration $q_i$ followed by a fibration $p_i$ using the factorizations of the model category \cat C. Then define $S(f) = \times_i p_i$. 

For every $i\in I$ there exists a pro-map $\varphi_i\colon \{X_i\}\rarrow X_i$ defined by the strict representative $\{\theta(i)=i,\; \{f_i = \id_{X_i}\colon X_i\rarrow X_i\}\}$. The same is true for $\{Y_i\}$. These maps define the canonical maps $\{X_i\}\rarrow \times_i X_i$ and $\{Y_i\}\rarrow \times_i Y_i$. Composition of the first map with $\times_i q_i$ finishes the definition of the morphism of maps $t_{\{f_i\}}\colon \{f_i\}\rarrow S(\{f_i\})$.

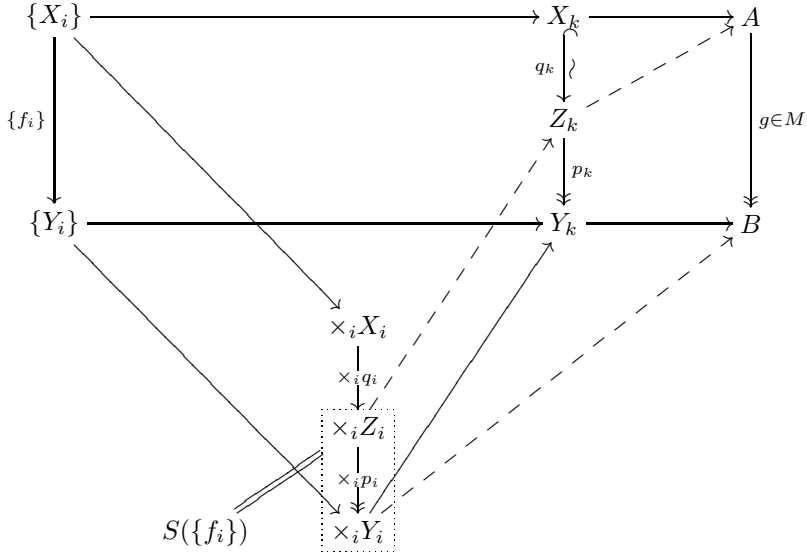
\[
 \xymatrix{
           \{X_i\} \ar[rrrr] \ar[dd]_{\{f_i\}}
                   \ar[dddrr]     & & & & X_k \ar[rr] \ar@{^{(}->}[d]^{\dir{~}}_{q_k}& & A \ar@{->>}[dd]^{g\in M}\\
                                  & & & & Z_k \ar@{->>}[d]^{p_k} \ar@{-->}[urr] \\
           \{Y_i\} \ar[rrrr]
             \ar[dddrr]          & & & & Y_k \ar[rr] & & B\\
             & & \times_i X_i \ar[d]|{\times_i q_i}
                    & &\\
             & & \times_i Z_i \ar@{->>}[d]|{\times_i p_i}
                    \ar@{-->}[uuurr] & &\\
             &S(\{f_i\})  & \times_i Y_i \ar[uuurr] \ar@{-->}[uuurrrr]& &
\save "5,3"."6,3"*[F.]\frm{} \ar@{=}"6,2" \restore
          }
\]

In order to verify the factorization property, fix an arbitrary map $\{f_i\}\rarrow g$ with $g\in M$ being a fibration between constant objects. It follows immediately from the definition of the morphism set between two pro-objects that any map into a constant pro-object factors through a map of the form $\varphi_i$ for some $i\in I$. Applying this to the category pro-$\Map(\cat C)\cong \Map(\text{pro-}\cat C)$, we find an index $k\in I$ such that the fixed map $\{f_i\}\rarrow g$ factors through $X_k\rarrow Y_k$.

In the diagram above the maps $\times_i Z_i\rarrow Z_k$ and $\times_i Y_i\rarrow Y_k$ are projections. The dashed map $\times_i Y_i\rarrow B$ is the composition of the projection with the map $Y_k\rarrow B$. Finally, the dashed map $Z_k\rarrow A$ is a lifting in the commutative square, which exists in the model category \cat C.

Now we are able to apply the non-functorial version of the generalized cosmall object argument to produce
for every map in pro-\cat C its factorization into a trivial cofibration followed by a fibration.

To obtain the second factorization we repeat the construction above for the class $N$ of trivial fibrations between constant objects and factorize all the maps $X_i\rarrow Y_i$ into cofibrations followed by trivial fibrations. Hence the generalized cosmall object argument may be applied to provide every map $f$ with a factorization into a cofibration followed by a trivial fibration in pro-\cat C.
\end{proof}

Let \cal S denote the category of simplicial sets with the standard model structure. Our next goal is to construct functorial factorizations in the localized model structure of \cite{Isaksen-sharp}. From now on the words \emph{cofibration}, \emph{fibration} and \emph{weak equivalence} refer to Isaksen's model structure.

Since the procedure of (left Bousfield) localization preserves the class of cofibrations and, hence, the class of trivial fibrations, the factorization into a cofibration followed by a trivial fibration was constructed in the theorem above. We keep the class $N$ of generating trivial fibrations the same as in the strict model structure. The class of generating fibrations $L$ is defined to be the class of all co-$n$-fibrations (see \cite[Def.~3.2]{Isaksen-sharp}) between constant pro-objects for all $n\in\NN$. 

\begin{proposition}
The class of trivial cofibrations equals $L$-{\rm proj}.
\end{proposition}
\begin{proof}
Every element of $L$ is a strong fibration \cite[Def.~6.5]{Isaksen-sharp}; therefore every trivial cofibration has the left lifting property with respect to $L$ by \cite[Prop.~14.5]{Isaksen-sharp}. Conversely, if a map $i$ has the left lifting property with respect to $L$, then $i$ has the left lifting property with respect to the class $L'$ of all retracts of $L$-cocell complexes. By \cite[Prop.5.2]{Isaksen-strict} $L'$ contains all strong fibrations. But then \cite[Prop.~6.6]{Isaksen-sharp} implies that  $L'$ contains all fibrations. Therefore, $i$ must be a trivial cofibration.
\end{proof}

\begin{theorem}\label{second}
Isaksen's model structure on {\rm pro-}\cal S is class-fibrantly generated with classes $L$ and $N$ of generating fibrations and generating trivial fibrations respectively.
\end{theorem}
\begin{proof}
It suffices to construct a factorization of every morphism of pro-\cal S into a trivial cofibration followed by a fibration. We apply the same construction as in Theorem~\ref{first} to the class $L$, except for the factorizations of the levelwise representation $\{f_i\}$. 

Apply first the \Mardesic functor of Theorem~\ref{Mardesic} in order to guarantee that our pro-system is indexed by a cofinite strongly directed set. Since the \Mardesic functor is naturally isomorphic to the identity, we abuse notation and keep calling the indexing category $I$. We construct the factorizations of the maps $f_i$ by induction on the number $n(i)$ of predecessors of $i$ and factor $f_i$ into an $n(i)$-cofibration $q_i$ followed by a co-$n(i)$-fibration, which is possible by \cite[Prop.~3.3]{Isaksen-sharp}. This is an ordinary induction on the set of natural numbers, since $I$ is now cofinite.

For any element $g\colon A\rarrow B$ of $L$, there is a number $n\in\NN$ such that $g$ is a co-$n$-fibration. We may always enlarge $k$ such that $n(k)\geq n$ and hence $q_k$ will be an $n$-cofibration by \cite[Lemma~3.6]{Isaksen-sharp}. Finally, the lift in the commutative square in the diagram in the proof of Theorem~\ref{first} exists by \cite[Def.~3.2]{Isaksen-sharp}.
\end{proof}

\bibliographystyle{abbrv}
\bibliography{xbib}

\end{document}